\documentclass{article}%
\usepackage{amssymb}
\usepackage{amsmath}
\usepackage{amsfonts}
\usepackage[dvips]{graphicx}%
\newtheorem{theorem}{Theorem}
\newtheorem{acknowledgement}[theorem]{Acknowledgement}

\newtheorem{proposition}[theorem]{Proposition}
\newtheorem{remark}[theorem]{Remark}

\newenvironment{proof}[1][Proof]{\noindent\textbf{#1.} }{\ \rule{0.5em}{0.5em}}
\begin{document}

\title{Smile Asymptotics II: Models with Known Moment Generating Function}
\author{Shalom Benaim and Peter Friz\\Statistical Laboratory, University of Cambridge}
\maketitle

\begin{abstract}
In a recent article the authors obtained a formula which relates explicitly
the tail of risk neutral returns with the wing behavior of the Black Scholes
implied volatility smile. In situations where precise tail asymptotics are
unknown but a moment generating function is available we first establish,
under easy-to-check conditions, tail asymptoics on logarithmic scale as soft
applications of standard Tauberian theorems. Such asymptotics are enough to
make the tail-wing formula work and we so obtain a version of Lee's moment
formula with the novel guarantee that there is indeed a limiting slope when
plotting implied variance against log-strike. We apply these results to
time-changed L\'{e}vy models and the Heston model. In particular, the
term-structure of the wings can be analytically understood.

\end{abstract}

\section{\bigskip Introduction}

Consider a random variable $X$ whose moment generating function (mgf)\ $M$ is
known in closed form, but whose density $f$ (if it exists) and distribution
function $F$ are, even asymptotically, unknown. For a large class of
distributions used for modelling (risk-neutral) returns in finance, $M$ is
finite only on part of the real line. Let us define $\bar{F}\equiv1-F $ and
$r^{\ast}$ as the least upper bound of all real $r$ for which $M\left(
r\right)  \equiv E[e^{rX}]<\infty$ and assume $r^{\ast}\in\left(
0,\infty\right)  $. An easy Chebyshev argument gives%
\begin{equation}
\lim\sup_{x\rightarrow\infty}\frac{-\log\bar{F}(x)}{x}=r^{\ast},
\label{limsupStmt}%
\end{equation}
but counter-examples show that the stronger statement%
\begin{equation}
-\log\bar{F}(x)\sim r^{\ast}x\text{ \ as }x\rightarrow\infty\text{ }
\label{as}%
\end{equation}
may not be true\footnote{We use the standard notation $g\left(  x\right)  \sim
h\left(  x\right)  \equiv\lim g\left(  x\right)  /h\left(  x\right)  =1$ as
$x\rightarrow\infty.$}. However, we do expect (\ref{as}) to be true if the
(right) tail of the distribution is reasonably behaved. Our interest in such
distributions stems from the fact that the crude tail asymptotics (\ref{as})
and the mild integrability condition $p^{\ast}=r^{\ast}-1>0$ are enough, via
the tail-wing formula \cite{BF}, to assert existence of a limiting slope of
Black Scholes implied variance $V^{2}$ as function of log-strike $k$. Indeed,
in standard notation, reviewed in section \ref{ApplSmile}, one has%
\begin{equation}
\lim_{k\rightarrow\infty}V^{2}(k)/k=2-4\left(  \sqrt{(p^{\ast})^{2}+p^{\ast}%
}-p^{\ast}\right)  . \label{rSlope}%
\end{equation}
Similarly, if $q^{\ast}\equiv\sup\left\{  q\in\mathbb{R}:M\left(  -q\right)
\equiv E[e^{-qX}]<\infty\right\}  \in(0,\infty)$ and the (left) tail is
reasonably behaved one expects $\log F(-x)\sim-q^{\ast}x$ as $x\rightarrow
\infty$ in which case the tail wing formula gives%
\begin{equation}
\lim_{k\rightarrow\infty}V^{2}\left(  -k\right)  /k=2-4\left(  \sqrt{(q^{\ast
})^{2}+q^{\ast}}-q^{\ast}\right)  . \label{lSlope}%
\end{equation}

It was already pointed out in \cite{BF} that the tail-wing formulae sharpen
Lee's celebrated moment formulae \cite{Lee, Ga}. In the present context, this
amounts to having a $\lim$ instead of a $\lim\sup$\footnote{Remark that, at
least when $p^{\ast}>0$, the moment formula is in fact recovered from the
tail-wing formula and (\ref{limsupStmt}).}. It must be noted that the
tail-wing formula requires some knowledge of the tails whereas the moment
formula is conveniently applicable by looking at the mgf (to obtain the
critical values $r^{\ast}$ and $-q^{\ast}$ ).

\bigskip

In this paper we develop criteria, checkable by looking \textit{a little
closer} at the mgf (near $r^{\ast}$ and $-q^{\ast}$), which will guarantee
that (\ref{rSlope}) resp. (\ref{lSlope}) hold. In view of the tail-wing
formula the problem is reduced to obtain criteria for (\ref{as}) resp. its
left-sided analogue. The proofs rely on Tauberian theorems and, as one
expects, the monograph \cite{BGT} is our splendid source.

The criteria are then fine-tuned to the fashionable class of time-changed
L\'{e}vy models \cite{S, CT} and checked explicitly for the examples of
Variance Gamma under Gamma-OU clock and Normal Inverse Gaussian with CIR
clock. We also check the criteria for the Heston model. In fact, it appears to
us that most (if not all) sensible models for stock returns with known mgf and
$p^{\ast},q^{\ast}\in\left(  0,\infty\right)  $ satisfy one of our criteria so
that (\ref{rSlope}) and (\ref{lSlope}) will hold.

Finally, we present some numerical results. The asymptotic regime becomes
visible for remarkably low log-strikes which underlines the practical value of
moment - and tail-wing formulae.

\section{Background\label{background} in Regular Variation}

\subsection{Asymptotic inversion}

If $f=f\left(  x\right)  $ is defined and locally bounded on $\,[X,\infty)$,
and tends to $\infty$ as $x\rightarrow\infty$ then the generalized inverse%
\[
f^{\leftarrow}\left(  x\right)  :=\inf\left\{  y\in\lbrack X,\infty):f\left(
y\right)  >x\right\}
\]
is defined on $[f\left(  X\right)  ,\infty)$ and is monotone increasing to
$\infty$. This applies in particular to $f\in R_{\alpha}$ with $\alpha>0$ and
Thm 1.5.12 in \cite{BGT} asserts that $f^{\leftarrow}\in R_{1/\alpha}$ and%
\[
f\left(  f^{\leftarrow}\left(  x\right)  \right)  \sim f^{\leftarrow}\left(
f\left(  x\right)  \right)  \sim x\text{ as }x\rightarrow\infty\text{.}%
\]
Given $f$ one can often compute $f^{\leftarrow}$ (up the asymptotic
equivalence) in terms of the \textit{Bruijn conjugate} of slowly varying
functions (Prop. 1.5.15, Section 5.2. and Appendix 5 in \cite{BGT}).

\subsection{Smooth Variation}

A positive function \thinspace$g$ defined in some neighbourhood of $\infty$
\textit{varies smoothly with index} $\alpha$, $g\in SR_{\alpha}$, iff
$h\left(  x\right)  :=\log\left(  g\left(  e^{x}\right)  \right)  $ is
$C^{\infty} $ and%
\[
h^{\prime}\left(  x\right)  \rightarrow\alpha,\,\,\,h^{\left(  n\right)
}\left(  x\right)  \rightarrow0\text{\thinspace\thinspace\thinspace for
}n=2,3,...\text{ as}\ x\rightarrow\infty.\text{ }%
\]

\begin{theorem}
[Smooth Variation Theorem, Thm 1.8.2 in \cite{BGT}]If $f\in R_{\alpha}$ then
there exist $f_{i}\in SR_{\alpha}$, $i=1,2$, with $f_{1}\sim f_{2}$ and
$f_{1}\leq f\leq f_{2}$ on some neighbourhood of $\infty$.
\end{theorem}

When $\alpha>0$ we can assume that $f_{1}$ and $f_{2}$ are strictly increasing
in some neighbourhood of $\infty$. In fact, we have

\begin{proposition}
\label{SVTcor}Let $\alpha>0$ and $g\in SR_{\alpha}$. Then $g$ is strictly
increasing in some neighbourhood of $\infty$ and $g^{\prime}\in SR_{\alpha-1}$.
\end{proposition}

\begin{proof}
By definition of $SR_{\alpha}$,%
\[
\frac{\partial}{\partial x}\log\left(  g\left(  e^{x}\right)  \right)
=\frac{1}{g\left(  e^{x}\right)  }g^{\prime}\left(  e^{x}\right)
e^{x}\rightarrow\alpha>0\text{ as }x\rightarrow\infty.
\]
This shows that, in some neighbourhood of $\infty$, $g^{\prime}$ is strictly
positive which implies that $g$ is strictly increasing. From Prop 1.8.1 in
\cite{BGT}, $g^{\prime}=\left\vert g^{\prime}\right\vert \in SR_{\alpha-1}$.
\end{proof}

\begin{remark}
In the situation of the last Proposition we have $\lim_{x\rightarrow\infty
}g\left(  x\right)  =\infty$ and hence, in some neighbourhood of $\infty$ ,
$g$ has a genuine inverse $g^{-1}$ which coincides with the generalized
inverse $g^{\leftarrow}$.
\end{remark}

\subsection{Exponential Tauberian Theory}

\begin{theorem}
[Kohlbecker's Theorem, Thm 4.12.1 and Cor 4.12.6 in \cite{BGT}]Let $U$ be a
non-decreasing right-continuous function on $\mathbb{R}$ with $U\left(
x\right)  =0$ for all $x<0$. Set%
\[
N\left(  \lambda\right)  :=\int_{[0,\infty)}e^{-x/\lambda}dU\left(  x\right)
,\,\,\,\lambda>0.
\]
Let $\alpha>1$ and $\chi\in R_{\alpha/\left(  \alpha-1\right)  }$. Then
\[
\log N\left(  \lambda\right)  \sim\left(  \alpha-1\right)  \chi\left(
\lambda\right)  /\lambda\text{ as }\lambda\rightarrow\infty
\]
iff%
\[
\log\mu\left[  0,x\right]  \sim\alpha x/\chi^{\leftarrow}\left(  x\right)
\text{ as }x\rightarrow\infty.
\]

\end{theorem}

\begin{theorem}
[Karamata's Tauberian Theorem, Thm 1.7.1 in \cite{BGT}]Let $U$ be a
non-decreasing right-continuous function on $\mathbb{R}$ with $U\left(
x\right)  =0$ for all $x<0$. If $l\in R_{0}$ and $c\geq0,\rho\geq0$, the
following are equivalent:%
\begin{align*}
U\left(  x\right)   &  \sim cx^{\rho}l\left(  x\right)  /\Gamma\left(
1+\rho\right)  \text{ as }x\rightarrow\infty\\
\hat{U}\left(  s\right)   &  \equiv\int_{0}^{\infty}e^{-sx}dU\left(  x\right)
\sim c s^{-\rho}l\left(  1/s\right)  \text{ as }s\rightarrow0+.
\end{align*}
(When $c=0$ the asymptotic relations are interpreted in the sense that
$U\left(  x\right)  =o\left(  x^{\rho}l\left(  x\right)  \right)  $ and
similar for $\hat{U}$.)
\end{theorem}

\begin{theorem}
[Bingham's Lemma, Thm 4.12.10 in \cite{BGT}]Let $f\in R_{\alpha}$ with
$\alpha>0$ such that that $e^{-f}$ is locally integrable at $+\infty$. Then%
\[
-\log\int_{x}^{\infty}e^{-f\left(  y\right)  }dy\sim f\left(  x\right)
\text{.}%
\]

\end{theorem}

\section{Moment generating functions and log-tails\label{MgfLT}}

Let $F$ be a finite Borel measure on $\mathbb{R}$, identified with its
(bounded, non-decreasing, right-continuous) distributions function, $F\left(
x\right)  \equiv F\left(  (-\infty,x]\right)  $. Its mgf is defined as%
\[
M\left(  s\right)  :=\int e^{sx}dF\left(  x\right)  .
\]
We define the critical exponents $q^{\ast}$ and $r^{\ast}$ via%
\[
-q^{\ast}\equiv\inf\left\{  s:M\left(  s\right)  <\infty\right\}  ,r^{\ast
}\equiv\sup\left\{  s:M\left(  s\right)  <\infty\right\}
\]
and make the \textbf{standing assumption} that%
\[
r^{\ast},q^{\ast}\in(0,\infty)\text{.}%
\]
In this section, we develop criteria which will imply the asymptotic
relations
\[
\log F\left(  (-\infty,-x]\right)  \sim-q^{\ast}x\text{, }\log F\left(
(x,\infty)\right)  \sim-r^{\ast}x\text{ as }x\rightarrow\infty.
\]

The assumption in the following Criterion I is simply that some derivative of
the mgf (at the critical exponent ) blows up in a regularly varying way.

\begin{theorem}
[Criterion I]Let $F$ be a bounded non-decreasing right-continuous function on
$\mathbb{R}$ and define $M=M\left(  s\right)  ,$ $q^{\ast}$ and $r^{\ast}$ as
above.\newline(i) If for some $n\geq0$, $M^{\left(  n\right)  }\left(
-q^{\ast}+s\right)  \sim s^{-\rho}l_{1}(1/s)$ for some $\rho>0$, $l_{1}\in
R_{0}$ as $s\rightarrow0+$ then%
\[
\log F\left(  (-\infty,-x]\right)  \sim-q^{\ast}x
\]
\newline(ii) If for some $n\geq0$, $M^{\left(  n\right)  }\left(  r^{\ast
}-s\right)  \sim s^{-\rho}l_{1}(1/s)$ for some $\rho>0$, $l_{1}\in R_{0}$ as
$s\rightarrow0+$ then%
\[
\log F\left(  (x,\infty)\right)  \sim-r^{\ast}x.
\]

\end{theorem}

\begin{proof}
Let us focus on case (ii), noting that case (i) is similar. We first discuss
$n=0$. The idea is an Escher-type change of measure followed by an application
of Karamata's Tauberian Theorem. We define a new measure $U$ on $[0,\infty)$
by a change-of-measure designed to get rid of the exponential decay,%
\[
dU\left(  x\right)  :=\exp\left(  r^{\ast}x\right)  dF\left(  x\right)  .
\]
We identify $U$ with its non-decreasing right-continuous distribution function
$x\mapsto U\left(  \left[  0,x\right]  \right)  $. The Laplace transform of
$U$ is given by%
\[
\hat{U}\left(  s\right)  =\int_{0}^{\infty}e^{-sx}dU\left(  x\right)
=\int_{0}^{\infty}e^{\left(  r^{\ast}-s\right)  x}dF\left(  x\right)
=M\left(  r^{\ast}-s\right)  -\int_{-\infty}^{0}e^{\left(  r^{\ast}-s\right)
x}dF\left(  x\right)
\]
so that%
\[
\left\vert \hat{U}\left(  s\right)  -M\left(  r^{\ast}-s\right)  \right\vert
\leq\int_{-\infty}^{0}e^{\left(  r^{\ast}-s\right)  x}dF\left(  x\right)  \leq
F\left(  0\right)  -F\left(  -\infty\right)  \leq2\left\Vert F\right\Vert
_{\infty}<\infty.
\]
Since $M\left(  r^{\ast}-s\right)  $ goes to $\infty$ as $s\rightarrow0+$ and
we see that $\hat{U}\left(  s\right)  \sim M\left(  r^{\ast}-s\right)  $ so
that $\hat{U}\in R_{\rho}$ as $s\rightarrow0$. Hence, there exists $l\in
R_{0}$ so that $\hat{U}\left(  s\right)  =\left(  1/s\right)  ^{\rho}l\left(
1/s\right)  $ and Karamata's Tauberian theorem tells us that $U\in R_{\rho}$,
namely%
\[
U\left(  x\right)  \sim x^{\rho}l\left(  x\right)  /\Gamma\left(
1+\rho\right)  \equiv x^{\rho}l^{\prime}\left(  x\right)  \text{ as
}x\rightarrow\infty
\]
where $l^{\prime}\in R_{0}$. Going back to the right-tail of $F$, we have for
$x\geq0$%
\[
F\left(  (x,\infty)\right)  =\int_{(x,\infty)}dF\left(  y\right)
=\int_{(x,\infty)}\exp\left(  -r^{\ast}y\right)  dU\left(  y\right)  \text{.}%
\]
We first assume that $U\in SR_{\rho}$. Under this assumption $U$ is smooth
with derivative $u=U^{\prime}$ $\in SR_{\rho-1}$ and we can write
\[
u\left(  y\right)  =y^{\rho-1}l^{\prime\prime}\left(  y\right)  \text{ with
}l^{\prime\prime}\in R_{0}.
\]
Then%
\begin{align*}
F\left(  (x,\infty)\right)   &  =\int_{(x,\infty)}\exp\left(  -r^{\ast
}y\right)  y^{\rho-1}l^{\prime\prime}\left(  y\right)  dy\\
&  =\int_{(x,\infty)}\exp\left[  -r^{\ast}y+\left(  \rho-1\right)  \log y+\log
l^{\prime\prime}\left(  y\right)  \right]  dy.
\end{align*}
Since $-\left[  -r^{\ast}y+\left(  \rho-1\right)  \log y+\log l^{\prime\prime
}\left(  y\right)  \right]  \sim r^{\ast}y\in R_{1}$ as $y\rightarrow\infty$
we can use Bingham's lemma to obtain%
\begin{equation}
-\log F\left(  (x,\infty)\right)  =-\log\int_{(x,\infty)}\exp\left[  -r^{\ast
}y+\left(  \rho-1\right)  \log y+\log l^{\prime\prime}\left(  y\right)
\right]  dU\left(  y\right)  \sim r^{\ast}y. \label{TailForSR}%
\end{equation}
We now deal with the general case of non-decreasing $U\in R_{\rho}$. From the
Smooth Variation Theorem and Proposition \ref{SVTcor} we can find $U_{-}%
,U_{+}\in SR_{\rho}$, strictly increasing in a neighbourhood of $\infty$, so
that%
\[
U_{-}\leq U\leq U_{+}\text{ and }U_{-}\sim U\sim U_{+}\text{.}%
\]
Below we use the change of variable $z=U\left(  y\right)  $ and $w=$
$U_{+}^{-1}\left(  z\right)  $. Noting that $U_{+}^{-1}\leq U^{\leftarrow}\leq
U_{-}^{-1}$ and using change-of-variable formulae, as found in \cite[p7-9]{RY}
for instance, we have
\begin{align*}
F\left(  (x,\infty)\right)   &  =\int_{(x,\infty)}\exp\left(  -r^{\ast
}y\right)  dU\left(  y\right) \\
&  =\int_{(U\left(  x\right)  ,\infty)}\exp\left(  -r^{\ast}U^{\leftarrow
}\left(  z\right)  \right)  dz\\
&  \leq\int_{(U\left(  x\right)  ,\infty)}\exp\left(  -r^{\ast}U_{+}%
^{-1}\left(  z\right)  \right)  dz\\
&  =\int_{(U_{+}^{-1}\left(  U\left(  x\right)  \right)  ,\infty)}\exp\left(
-r^{\ast}w\right)  dU_{+}\left(  w\right)  \text{.}%
\end{align*}
Similar to the derivation of (\ref{TailForSR}), Bingham's lemma leads to%
\[
-\log\int_{(U_{+}^{-1}\left(  U\left(  x\right)  \right)  ,\infty)}\exp\left(
-r^{\ast}w\right)  dU_{+}\left(  w\right)  \sim r^{\ast}U_{+}^{-1}\left(
U\left(  x\right)  \right)  .
\]
Noting that $U_{+}^{-1}$ is non-decreasing, $U_{+}^{-1}\left(  U\left(
x\right)  \right)  \leq U_{+}^{-1}\left(  U_{+}\left(  x\right)  \right)  =x$
so that\footnote{By $g\lesssim h$ we mean $\lim\sup f\left(  x\right)
/g\left(  x\right)  \leq1$ as $x\rightarrow\infty$.}%
\[
-\log F\left(  [x,\infty)\right)  \lesssim r^{\ast}x
\]
The same argument gives the lower bound $-\log F\left(  (x,\infty)\right)
\gtrsim r^{\ast}x$ and we conclude that%
\[
-\log F\left(  (x,\infty)\right)  \sim r^{\ast}x.
\]
We now show how $n>0$ follows from $n=0$. Define $V$ on $[0,\infty)$ by
\[
dV\left(  x\right)  :=x^{n}dF\left(  x\right)  .
\]
Clearly, $V$ induces a non-decreasing, right continuous distribution on
$\mathbb{R}$, $V\left(  x\right)  :=V\left(  [0,x]\right)  $ for $x\geq0$ and
$V\left(  x\right)  \equiv0$ for $x<0$. The distribution function $V\left(
x\right)  $ is also bounded since%
\[
\int_{0}^{\infty}x^{n}dF\left(  x\right)  <\infty
\]
which follows a forteriori from the standing assumption of exponential
moments. We will write $\bar{V}(x)$ for $V(x,\infty).$

Note that $V$ has a mgf $M_{V}\left(  s\right)  $, finite at least for
$s\in(0,r^{\ast})$, given by%
\begin{align*}
M_{V}\left(  s\right)   &  \equiv\int e^{sx}dV\left(  x\right)  =\int
_{0}^{\infty}x^{n}e^{sx}dF\\
&  =\int x^{n}e^{sx}dF+C=M^{(n)}\left(  s\right)  +C
\end{align*}
where\footnote{One could do without the assumption $\int_{-\infty}%
^{0}\left\vert x\right\vert dF$ (which follows a forteriori from the standing
assumption $q^{\ast}>0$). Finiteness of $F$ on $\left(  -\infty,0\right)  $ is
enough.}
\[
0\leq C\equiv-\int_{-\infty}^{0}x^{n}e^{sx}dF\leq\int_{-\infty}^{0}\left\vert
x\right\vert ^{n}dF<\infty.
\]
By assumption, $M^{(n)}$ is regularly varying with index $\rho$ at $r^{\ast}$
and it follows that, as $s\rightarrow0+$,%
\[
M_{V}\left(  r^{\ast}-s\right)  =M^{(n)}\left(  r^{\ast}-s\right)  +O\left(
1\right)  \sim s^{-\rho}l_{1}(1/s).
\]
We now use the "$n=0$" result on the distribution function $V$ resp. its mgf
$M_{V}$ and obtain%
\[
-\log V\left(  [x,\infty)\right)  \equiv-\log\bar{V}\left(  x\right)  \sim
r^{\ast}x\in R_{1}%
\]
Assume first that $-\log\bar{V}\left(  x\right)  \in SR_{1}$. Then $V$ has a
density $V^{\prime}\equiv v$ and%
\[
v\left(  x\right)  =\partial_{x}(V(\infty)-\bar{V}(x))=-\bar{V}\left(
x\right)  \partial_{x}\left(  \log\bar{V}\left(  x\right)  \right)  \sim
r^{\ast}\bar{V}\left(  x\right)  \text{ as }x\rightarrow\infty
\]
since functions in $SR_{1}$ are stable under differentiation in the sense that
$\partial_{x}\left(  -\log\bar{V}\left(  x\right)  \right)  \sim\partial
_{x}\left(  r^{\ast}x\right)  =r^{\ast}$. In particular, we have $\log
v\left(  x\right)  \sim\log\bar{V}\left(  x\right)  \sim-r^{\ast}x$. After
these preparations we can write
\begin{align*}
F\left(  (x,\infty)\right)   &  =\int_{(x,\infty)}dF\left(  y\right) \\
&  =\int_{(x,\infty)}\frac{1}{y^{n}}v\left(  y\right)  dy\\
&  =\int_{(x,\infty)}\exp\left[  \log v\left(  y\right)  -n\log y\right]  dy
\end{align*}
and Bingham's lemma implies that $\log F\left(  (x,\infty)\right)
\sim-r^{\ast}x$. The general case of $\log\bar{V}\left(  x\right)  \in R_{1}$
follows by a smooth variation and comparison argument as earlier.
\end{proof}

The next criterion deals with exponential blow-up of $M$ at its critical values.

\begin{theorem}
[Criterion II]Let $F,M,q^{\ast},r^{\ast}$ be as above.\newline(i) If $\log
M\left(  -q^{\ast}+s\right)  \sim s^{-\rho}l_{1}(1/s)$ for some $\rho>0$,
$l_{1}\in R_{0}$ as $s\rightarrow0+$ then
\[
\log F\left(  (-\infty,-x]\right)  \sim-q^{\ast}x
\]
\newline(ii) If $\log M\left(  r^{\ast}-s\right)  \sim s^{-\rho}l_{1}(1/s)$
for some $\rho>0$, $l_{1}\in R_{0}$ as $s\rightarrow0+$ then%
\[
\log F\left(  (x,\infty)\right)  \sim-r^{\ast}x.
\]

\end{theorem}

\begin{proof}
As for Criterion I, the idea is an Escher-type change of measure followed by a
suitable Tauberian theorem; in the present case we need Kohlbecker's Theorem.
Let us focus on case (ii), noting that case (i) is similar. A new measure $U$
on $[0,\infty)$ is defined by
\[
dU\left(  x\right)  :=\exp\left(  r^{\ast}x\right)  dF\left(  x\right)  .
\]
We identify $U$ with its non-decreasing right-continuous distribution function
$x\mapsto U\left(  \left[  0,x\right]  \right)  $ and define the transform%
\[
N\left(  \lambda\right)  =\int_{0}^{\infty}e^{-x/\lambda}dU\left(  x\right)
=\int_{0}^{\infty}e^{\left(  r^{\ast}-1/\lambda\right)  x}dF\left(  x\right)
=M\left(  r^{\ast}-1/\lambda\right)  -\int_{-\infty}^{0}e^{\left(  r^{\ast
}-1/\lambda\right)  x}dF\left(  x\right)
\]
so that%
\[
\left\vert N\left(  \lambda\right)  -M\left(  r^{\ast}-1/\lambda\right)
\right\vert \leq\int_{-\infty}^{0}e^{\left(  r^{\ast}-1/\lambda\right)
x}dF\left(  x\right)  \leq F\left(  0\right)  -F\left(  -\infty\right)
\leq2\left\Vert F\right\Vert _{\infty}<\infty.
\]
Thus,%
\[
N\left(  \lambda\right)  =M\left(  r^{\ast}-1/\lambda\right)  +O\left(
1\right)  \text{ as }\lambda\rightarrow\infty
\]
and, in particular, since $\lim_{\lambda\rightarrow\infty}\log M\left(
r^{\ast}-1/\lambda\right)  =\lim_{\lambda\rightarrow\infty}M\left(  r^{\ast
}-1/\lambda\right)  =\infty$ from the assumption (ii) we see that
\[
\log N\left(  \lambda\right)  \sim\log M\left(  r^{\ast}-1/\lambda\right)
\text{ }\sim\lambda^{\rho}l_{1}(\lambda)\in R_{\rho}\text{ as }\lambda
\rightarrow\infty.
\]
Define $\alpha\in(1,\infty)$ as the unique solution to $\rho+1=\alpha/\left(
\alpha-1\right)  $ and note%
\[
\chi\left(  \lambda\right)  :=\frac{\lambda}{\left(  \alpha-1\right)  }\log
N\left(  \lambda\right)  \in R_{\rho+1}=R_{\alpha/\left(  \alpha-1\right)  }.
\]
Using that $\chi^{\leftarrow}\in R_{\left(  \alpha-1\right)  /\alpha
}=R_{1-1/\alpha}$, Kohlbecker's Tauberian Theorem tells us that%
\[
\log U\left(  \left[  0,x\right]  \right)  \equiv\log U\left(  x\right)
\sim\alpha x/\chi^{\leftarrow}\left(  x\right)  \in R_{1/\alpha}\text{ as
}x\rightarrow\infty.
\]
In particular, there exists $l\in R_{0}$ so that $\log U\left(  x\right)
=\alpha x^{1/\alpha}l\left(  x\right)  $. We first assume that $\log U\in
SR_{1/\alpha}$. Then $U$ has a density $u\left(  .\right)  \in SR_{1/\alpha
-1}$ and%
\[
u\left(  x\right)  =U\left(  x\right)  \partial_{x}\left(  \log U\left(
x\right)  \right)  \sim U\left(  x\right)  x^{1/\alpha-1}l\left(  x\right)  .
\]
In particular,%
\[
\log u\left(  x\right)  \sim\log U\left(  x\right)  \in R_{1/\alpha}\text{ as
}x\rightarrow\infty.
\]
Now, $y\mapsto r^{\ast}y\in R_{1}$ dominates $R_{1/\alpha}$ (since
$1/\alpha<1$) in the sense that%
\[
r^{\ast}y-\log u\left(  y\right)  \sim r^{\ast}y.
\]
Thus, from%
\begin{align*}
F\left(  (x,\infty)\right)   &  =\int_{(x,\infty)}dF\left(  y\right)
=\int_{[x,\infty)}\exp\left(  -r^{\ast}y\right)  u\left(  y\right)  dy\\
&  =\int_{(x,\infty)}\exp\left[  -r^{\ast}y+\log u\left(  y\right)  \right]
\end{align*}
and Bingham's lemma we deduce that%
\[
-\log F\left(  (x,\infty)\right)  \sim r^{\ast}x\text{.}%
\]
The general case, $\log U\in R_{1/\alpha}$, is handled via smooth variation as
earlier. Namely, we can find smooth minorizing and majorizing functions for
$\log U$, say $G\_$ and $G_{+}$. After defining $U_{\pm}=\exp G_{\pm}$ we have%
\[
\log U_{-}\sim\log U\sim\log U_{+}\text{ and }U_{-}\leq U\leq U_{+}.
\]
Then, exactly as in the last step of the proof of Criterion I,%
\[
F\left(  (x,\infty)\right)  =\int_{(x,\infty)}\exp\left(  -r^{\ast}y\right)
dU\left(  y\right)  \leq\int_{(U_{+}^{-1}\left(  U\left(  x\right)  \right)
,\infty)}\exp\left(  -r^{\ast}w\right)  dU_{+}\left(  w\right)
\]
and from Bingham's lemma,%
\[
-\log F\left(  (x,\infty)\right)  \lesssim r^{\ast}U_{+}^{-1}\left(  U\left(
x\right)  \right)  \sim r^{\ast}x\text{.}%
\]
Similarly, $-\log F\left(  (x,\infty)\right)  \gtrsim r^{\ast}x$ and the proof
is finished.
\end{proof}

\section{\bigskip Application to Smile Asymptotics\label{ApplSmile}}

We start with a few recalls to settle the notation. The normalized price of a
Black-Scholes call with log-strike $k$ is given by%
\[
c_{BS}\left(  k,\sigma\right)  =\Phi\left(  d_{1}\right)  -e^{k}\Phi\left(
d_{2}\right)  \,
\]
with $d_{1,2}\left(  k\right)  =-k/\sigma\pm\sigma/2$. If one models
risk-neutral returns with a distribution function $F$, the implied volatility
is the (unique) value $V\left(  k\right)  $ so that
\[
c_{BS}\left(  k,V\left(  k\right)  \right)  =\int_{k}^{\infty}\left(
e^{x}-e^{k}\right)  dF\left(  x\right)  =:c\left(  k\right)  \text{.}%
\]
Set $\psi\left[  x\right]  \equiv2-4\left(  \sqrt{x^{2}+x}-x\right)  $ and
recall $\bar{F}\equiv1-F$. The following is a special case of the tail-wing
formula \cite{BF}.

\begin{theorem}
\label{FromBF}\bigskip Assume that $-\log F\left(  -k\right)  /k\sim q^{\ast}$
for some $q^{\ast}\in\left(  0,\infty\right)  .$ Then%
\[
V(-k)^{2}/k\sim\psi\left[  -\log F\left(  -k\right)  /k\right]  \sim
\psi\left(  q^{\ast}\right)  .
\]
Similarly, assume that $-\log\bar{F}\left(  k\right)  /k\sim p^{\ast}+1$ for
some $p^{\ast}\in\left(  0,\infty\right)  $. Then%
\[
V(k)^{2}/k\sim\psi\left[  -1-\log\bar{F}\left(  k\right)  /k\right]  \sim
\psi\left(  p^{\ast}\right)  .
\]

\end{theorem}

As earlier, let $M\left(  s\right)  =\int\exp\left(  sx\right)  dF\left(
x\right)  $ denote the mgf of risk-neutral returns and now \textit{define} the
critical exponents $r^{\ast}$and $-q^{\ast}$ exactly as in the beginning of
the last section \ref{MgfLT}. Combining the results therein with Theorem
\ref{FromBF} we obtain

\begin{theorem}
\label{main}\bigskip If $q^{\ast}\in\left(  0,\infty\right)  $ and $M$
satisfies part (i) of Criteria I or II then%
\[
V(-k)^{2}/k\sim\psi\left(  q^{\ast}\right)  \text{ as }k\rightarrow
\infty\text{.}%
\]
Similarly, if $\ r^{\ast}\equiv p^{\ast}+1\in\left(  1,\infty\right)  $ and
$M$ satisfies part (ii) of Criteria I or II then%
\[
V(k)^{2}/k\sim\psi\left(  p^{\ast}\right)  \text{ as }k\rightarrow
\infty\text{.}%
\]

\end{theorem}

\section{First \bigskip Examples\label{Examples}}

The examples discussed in this section model risk-neutral log-price by
L\'{e}vy processes and there is no loss of generality to focus on unit
time.\footnote{In fact, L\'{e}vy models that satisfy one of our criteria have
no term structure of implied variance slopes.}

\subsection{Criterion I with n=0: the Variance Gamma Model}

The Variance Gamma model $VG=VG\left(  m,g,C\right)  $ has mgf%
\[
M(s)=\left(  \frac{gm}{gm+(m-g)s-s^{2}}\right)  ^{C}=\left(  \frac
{gm}{(m-s)(s+g)}\right)  ^{C}\text{.}%
\]
The critical exponents are obviously given by $r^{\ast}=m$ and
$q^{\ast}=g$.
Focusing on the first, we have%
\[
M(r^{\ast}-s)\sim\left(  \frac{gm}{m+g}\right)  ^{C}s^{-C}\text{ as
}s\rightarrow0+
\]
which shows the Criterion I is satisfied with $n=0$. Theorem \ref{main} now
identifies the asymptotic slope of the implied variance to be $\psi\left(
r^{\ast}-1\right)  =\psi\left(  m-1\right)  $. Similarly, the left slope is
seen to be $\psi\left(  q^{\ast}\right)  =\psi\left(  g\right)  $. We remark
that \cite{AB1} contains tail estimates for $VG$ which lead, via the tail-wing
formula, to the same result.

\subsection{\bigskip Criterion\ I with n%
$>$%
0: the Normal Inverse Gaussian Model}

The Normal Inverse Gaussian Model $NIG=$ $NIG\left(  \alpha,\beta,\mu
,\delta\right)  $ has mgf given by%
\[
M\left(  s\right)  =\exp\left\{  \delta\left\{  \sqrt{\alpha^{2}-\beta^{2}%
}-\sqrt{\alpha^{2}-\left(  \beta+s\right)  ^{2}}\right\}  +\mu s\right\}  .
\]
By looking at the endpoints of the strip of analyticity the critical exponents
are immediately seen to be $r^{\ast}=$ $\alpha-\beta,\,q^{\ast}=\alpha+\beta$
and we focus again on the first. While $M(s)$ converges to the finite constant
$M(r^{\ast})$ as $s\rightarrow r^{\ast}-$ we have%
\begin{align*}
M^{\prime}(s)/M\left(  s\right)   &  =(2\delta(\beta+s)[\alpha^{2}%
-(\beta+s)^{2}]^{-1/2}+\mu)\\
\text{and }M^{\prime}(r^{\ast}-s)  &  \sim2\delta\alpha\sqrt{2\alpha}%
s^{-1/2}M(r^{\ast})\text{ as }s\rightarrow0+\text{.}%
\end{align*}

We see that Criterion I is satisfied with $n=1$ and Theorem \ref{main} gives
the asymptotic slope $\psi\left(  r^{\ast}-1\right)  =\psi\left(  \alpha
-\beta-1\right)  $. Similarly, the left slope is seen to be $\psi\left(
q^{\ast}\right)  =\psi\left(  \alpha+\beta\right)  $. We remark that the same
slopes were computed in \cite{BF} using the tail-wing formula and explicitly
known density asymptotics for $NIG$.

\subsection{Criterion\ II: the Double Exponential Model}

The double exponential model $DE=DE\left(  \sigma,\mu,\lambda,p,q,\eta
_{1},\eta_{2}\right)  $ has mgf%
\[
\log M(s)=\frac{1}{2}\sigma^{2}s^{2}+\mu s+\lambda\left(  \frac{p\eta_{1}%
}{\eta_{1}-s}+\frac{q\eta_{2}}{\eta_{2}+s}-1\right)  .
\]
Clearly, $r^{\ast}=\eta_{1}$ and as $s\rightarrow0+$%
\[
\log M(\eta_{1}-s)\sim\frac{1}{2}\sigma\eta_{1}^{2}+\mu\eta_{1}+\lambda\left(
\frac{p\eta_{1}}{s}+\frac{q\eta_{2}}{\eta_{2}+\eta_{1}}-1\right)  \sim\lambda
p\eta_{1}s^{-1}%
\]
and we see that Criterion II is satisfied. As above, this implies asympotic
slopes $\psi\left(  r^{\ast}-1\right)  =\psi\left(  \eta_{1}-1\right)  $ on
the right and $\psi\left(  \eta_{2}\right)  $ on the left.

\section{\bigskip Time changed L\'{e}vy processes}

We now discuss how to apply our results to time changed L\'{e}vy processes
\cite{S, SST, CT}. \ To do this, we only need to check that the moment
generating function of the marginals of the process, satisfies one of the
three criteria.

To this end, consider a L\'{e}vy process $L=L\left(  t\right)  $ described
through its cumulant generating function (cgf) $K_{L}$ at time $1$, that is,%
\[
K_{L}\left(  v\right)  =\log E\left[  \exp\left(  vL_{1}\right)  \right]
\]
and an independent random clock $T=T\left(  \omega\right)  \geq0$ with cgf
$K_{T}$. It follows that the mgf of $L\circ T$ is given by%
\[
M(v)=\mathbb{E}\left[  \mathbb{E}\left(  e^{vL_{T}}|T\right)  \right]
=\mathbb{E}\left[  e^{K_{L}\left(  v\right)  T}\right]  =\exp\left[
K_{T}(K_{L}(v)\right]  .
\]
Therefore, in order to apply our Theorem \ref{main} to time-changed L\'{e}vy
models we need to check if $M=\exp\left[  K_{T}(K_{L}(\cdot)\right]  $
satisfies criterion I or II so that $-\log\bar{F}\left(  x\right)  /x$ tends
to a positive constant. Here, as earlier, $F$ denotes the distribution
function of $M$ and $\bar{F}\equiv1-F$. The following theorem gives sufficient
conditions for this in terms of $K_{T}$ and $K_{L}$. We shall write
$M_{T}\equiv\exp\left(  K_{T}\right)  $ and $M_{L}\equiv\exp\left(
K_{L}\right)  $. For brevity, we only discuss the right tail\footnote{In fact,
the elegant change-of-measure argument in Lee \cite{Lee} allows a formal
reduction of the left tail behaviour to the right tail behaviour.} and set%
\[
p_{L}=\sup\left\{  s:M_{L}\left(  s\right)  <\infty\right\}  ,\,\,\,p_{T}%
=\sup\left\{  s:M_{T}\left(  s\right)  <\infty\right\}  .
\]

\begin{theorem}
\label{TCLT}With notation as above, assuming $p_{L},p_{T}>0$, we
have:\newline(i.1) If $K_{L}(p)=p_{T}$ for some $p\in\lbrack0,p_{L})$ and
$M_{T}$ satisfies either criterion then%
\[
\log\bar{F}(x)\sim-px.
\]
(i.2) If $K_{L}(p)=p_{T}$ for $p=p_{L}$ and $M_{T},M_{L}$ satisfy
either
criterion then%
\[
\log\bar{F}(x)\sim-px.
\]
(ii) If $K_{L}(p)<p_{T}$ for all $p\in\lbrack0,p_{L}]$ and $M_{L}$ satisfies
either criterion then%
\[
\log\bar{F}(x)\sim-p_{L}x.
\]

\end{theorem}

\begin{remark}
\label{UniqueSolForP}It is worth noting that there cannot be more
than one solution to $K_{L}(p)=p_{T}$. To see this, take any $v$
such that $v>0$ and $K_{L}\left(  v\right)  >0$ (any solution to
$K_{L}(p)=p_{T}>0$ will satisfy this!) From $M_{L}\equiv\exp
K_{L}$ it follows that $M_{L}\left(  0\right) =1$ and $M_{L}\left(
v\right)  >1$. By convexity of $M_{L}\left( \cdot\right)  $ it is
easy to see that $M_{L}^{\prime}\left(  v\right)  $ is strictly
positive and the same is true for $K_{L}^{\prime}\left(  v\right)
=M_{L}^{\prime}\left(  v\right)  /M_{L}\left(  v\right)  $. It
follows that $w\geq v\implies$ $K_{L}\left(  w\right)  \geq
K_{L}\left(  v\right)  >0$ and the set of all $v>0:K_{L}\left(
v\right)  >0$ is connected and $K_{L}$ restricted to this set is
strictly increasing. This shows that there is at most one solution
to $K_{L}(p)=p_{T}$.
\end{remark}

\begin{proof}
(i.1) Noting that $p>0$ let as first assume that $M_{T}$ satisfies criterion I
(at $K_{L}(p)=p_{T}$ with some $n\geq0\,$) so that for some $\rho>0$ and $l\in
R_{0},$%
\[
M_{T}^{\left(  n\right)  }\left(  u\right)  \sim\left(  p_{T}-u\right)
^{-\rho}l\left(  \left(  p_{T}-u\right)  ^{-1}\right)  \text{ as }u\uparrow
p_{T}.
\]
From $M=M_{T}\circ K_{L}$ we have $M^{\prime}=M_{T}^{\prime}\left(
K_{L}\right)  K_{L}^{\prime}$ and, by iteration, $M^{\left(  n\right)  }$
equals $M_{T}^{\left(  n\right)  }\left(  K_{L}\right)  \left(  K_{L}^{\prime
}\right)  ^{n}$ plus a polynomial in $M_{T}\left(  \cdot\right)
,...,M_{T}^{\left(  n-1\right)  }\left(  \cdot\right)  $ which remains bounded
when the argument approaches $p_{T}$. Noting that $K_{L}^{\prime}\left(
p\right)  >0$ (see remark above) we absorb the factor $[K_{L}^{\prime}\left(
p\right)  ]^{n}$ into the slowly varying function and see that
\[
M^{\left(  n\right)  }(v)\sim(p_{T}-K_{L}(v))^{-\rho}l((p_{T}-K_{L}%
(v))^{-1})\text{ for }\rho\text{ as above and some }l\text{ }\in R_{0}%
\]
as $K_{L}\left(  v\right)  $ tends to $p_{T}$ which follows from $v\uparrow
p$. Using analyticity of $K_{L}$ in $(0,p_{L})$ and $K_{L}^{\prime}$ $\left(
p\right)  \neq0$ it is clear that $p_{T}-K_{L}(v)\sim K_{L}^{\prime}(p)(p-v)$
as $v\uparrow p$ and so%
\[
M^{\left(  n\right)  }(p-v)\sim K_{L}^{\prime}(p)^{-\rho}v^{-\rho}l(1/v)\text{
as }v\rightarrow0+.
\]
This shows that $M$ satisfies Criterion I (with the same $n$ as $M_{T}$). A
similar argument shows that $M$ satisfies Criterion II if $M_{T}$ does. Either
way, the assert tail behaviour of $\log\bar{F}$ follows.\newline(i.2) The
(unlikely!)\ case $K_{L}(p_{L})=p_{T}$ involves similar ideas and is left to
the reader.\newline(ii)\ We now assume that%
\[
\sup_{p\in\left[  0,p_{L}\right]  }K_{L}(p)<p_{T}<\infty.
\]
and $M_{L}$ satisfies either criterion (at $p_{L}$). Since $M_{L}=\exp K_{L}$
stays bounded as its argument approaches the critical value $p_{L}$ it is
clear that $M_{L}$ cannot satisfy criterion II or criterion I with $n=0$ and
there must exist a smallest integer $n$ such that%
\[
M_{L}^{(n)}(p_{L}-x)\sim x^{-\rho}l(x)\text{ \ \ as }x\rightarrow p_{L}%
\]
for some $\rho>0$ and $l\in R_{0}$. We note that%
\[
M^{(n)}(v)=(K_{L}^{(n)}(v)K_{T}^{\prime}(K_{L}(v))+f(v))\exp(K_{T}(K_{L}(v)))
\]
where $f(v)$ is a polynomial function of the first $\left(  n-1\right)  $
derivatives of $K_{L}$ and the first $n$ derivatives of $K_{T}$ evaluated at
$K_{L}(v)$, which are all bounded for $0\leq v\leq p_{L}$. Noting that
positivity of $T$ implies $M_{T}^{\prime}>0$ and hence $K_{T}^{\prime}>0$ we
see that as $v\uparrow p_{L}$%
\[
M^{(n)}(v)\sim K_{L}^{(n)}(v)K_{T}^{\prime}(K_{L}(p_{L}))M\left(
p_{L}\right)  .
\]
Applying this to $K_{T}\left(  x\right)  \equiv x$ leads immediately to%
\[
K_{L}^{(n)}(v)\sim M_{L}^{(n)}(v)/M_{L}(v)\sim x^{-\rho}l(x)/M_{L}(p_{L}).
\]
as $v\uparrow p_{L}$, and so $M$ satisfies criterion I.\newline
\end{proof}

We now discuss examples to which the above analysis is applicable. For all
examples we plot the total variance smile\footnote{That is, $V^{2}\left(
k,t\right)  \equiv\sigma^{2}\left(  k,t\right)  t.$} for several maturities
and compare with straight lines\footnote{These lines have been
parallel-shifted so that they are easier to compare with the actual smile.} of
correct slope as predicted by Theorem \ref{main}. All plots are based on the
calibrations obtained in \cite{SST}. This is also where the reader can find
more details about the respective model parameters.

\subsection{Variance Gamma with Gamma-OU time change}

We will consider the Variance Gamma process with a Gamma-Ornstein-Uhlenbeck
time change and refer to \cite{SST} for details. From earlier, the Variance
Gamma process has cumulant generating function%
\[
K_{L}(v)=C\log\left(  \frac{gm}{(m-v)(v+g)}\right)  \text{ for }v\in(-g,m)
\]
We note that $K_{L}\left(  \left[  0,m\right]  \right)  =\left[
0,\infty\right]  $ so that $p_{L}=\infty$. The Gamma-Ornstein-Uhlenbeck clock
$T=T\left(  \omega,t\right)  $ has cgf%
\[
K_{T}(v)=vy_{0}\lambda^{-1}(1-e^{-\lambda t})+\frac{\lambda a}{v-\lambda
b}\left[  b\log\left(  \frac{b}{b-v\lambda^{-1}(1-e^{-\lambda t})}\right)
-vt\right]
\]
We need to examine how this function behaves around the endpoint of its strip
of regularity. At first glance, it appears that the function tends to infinity
as $v\uparrow\lambda b$, because of the $\frac{\lambda a}{v-\lambda b}$ term.
However, upon closer examination, we can see that this is in fact a removable
singularity, and the term of interest to us is the $\log(...)$ term. This term
tends to infinity as $v\rightarrow\lambda b(1-e^{-\lambda t})^{-1}=:p_{T}$.
After some simple algebra, we see that%
\begin{align*}
e^{K_{T}(v)} &  =\left(  \frac{b}{b-v\lambda^{-1}(1-e^{-\lambda t})}\right)
^{\frac{\lambda ab}{v-\lambda b}}\exp\left\{  vy_{0}\lambda^{-1}(1-e^{-\lambda
t})-\frac{vt\lambda a}{v-\lambda b}\right\}  \\
&  \sim\left(  \frac{p_{T}}{p_{T}-v}\right)  ^{\frac{\lambda ab}{p_{T}-\lambda
b}}\exp\left\{  p_{T}y_{0}\lambda^{-1}(1-e^{-\lambda t})-\frac{p_{T}t\lambda
a}{p_{T}-\lambda b}\right\}  \text{ as }v\uparrow p_{T}.
\end{align*}
Therefore, $\exp\left(  K_{T}\right)  $ satisfies Criterion I with $n=0$ and
part (i.1) of Theorem \ref{TCLT} shows that $M$ does too and so that $\log
\bar{F}(x)\sim-px$ where $p$ is determined by the equation%
\[
K_{L}(p)=p_{T}=\lambda b(1-e^{-\lambda t})^{-1}%
\]
and can be calculated explicitly,%
\[
p=\frac{m-g+\sqrt{(m-g)^{2}+4gm(1-\exp(-\lambda b/C(1-e^{-\lambda t}))}}{2}.
\]

\begin{figure}[tbh]
\includegraphics{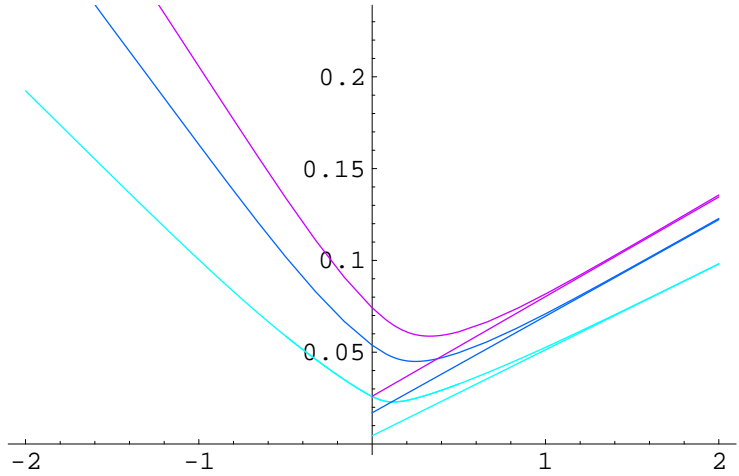}\caption{VG with Gamma-OU time change. Parameters from
\cite{SST}. Total implied variance and slopes for three maturities $t=0.4,0.9$
and $1.3$ years.}%
\end{figure}

\subsection{\bigskip Normal Inverse Gaussian with CIR time change}

The cgf of the Cox-Ingersoll-Ross (CIR)\ clock $T=T\left(  \omega,t\right)  $
is given by
\[
K_{T}(v)=\kappa^{2}\eta t/\lambda^{2}+2y_{0}v/(\kappa+\gamma\coth(\gamma
t/2))-\frac{2\kappa\eta}{\lambda^{2}}\log[\sinh\gamma t/2(\coth\gamma
t/2+\frac{\kappa}{\gamma})]\text{ }%
\]
where%
\[
\gamma=\sqrt{\kappa^{2}-2\lambda^{2}v}\text{.}%
\]
This clearly tends to infinity as $I(v)\equiv\kappa+\gamma(v)\coth
(\gamma(v)t/2)\rightarrow0$, and we can define $p_{T}$ as solution to the
equation $I\left(  p_{T}\right)  =0$. Using l'H\^{o}pital's rule, it is easy
to check that%
\[
\frac{p_{T}-v}{\kappa+\gamma(v)\coth(\gamma(v)t/2)}t
\]
tends to a constant as $v\rightarrow p_{T}$ , and so $2y_{0}v/(\kappa
+\gamma\coth(\gamma t/2)$ is regularly varying of index $1$ as a\ function of
$\left(  p_{T}-v\right)  ^{-1}$. It is clear that this is the dominant term in
this limit, and so $M_{T}\equiv\exp\left(  K_{T}\right)  $ satisfies criterion
II (at $p_{T}$). From earlier, the NIG cgf is\footnote{Following \cite{SST} we
take $\mu=0$ here.}
\[
K_{L}(v)=-\delta(\sqrt{\alpha^{2}-(\beta+v)^{2}}-\sqrt{\alpha^{2}-\beta^{2}%
})\text{ for }v\leq\alpha-\beta
\]
from which we see that $p_{L}=\alpha-\beta>0$ and%

\[
\sup_{v\in\lbrack0,\alpha-\beta]}K_{L}(v)=\delta\sqrt{\alpha^{2}-\beta^{2}%
}\text{.}%
\]
Therefore, the behavior of $M$ on the edge of the strip of analyticity, and
the location of the critical value, will depend on whether this supremum is
more or less than $p_{T}$; if it is less than $p_{T}$, the latter is never
reached. Recalling that $\exp\left(  K_{L}\right)  $ satisfies Criterion I
with $n=1$, we apply part (ii) of Theorem \ref{TCLT} and obtain%
\[
-\log\bar{F}(x)\sim p_{L}x=(\alpha-\beta)x.
\]
Otherwise, there exists $p\in(0,\alpha-\beta]$ such that $K_{L}(p)=p_{T}$, for
some $p\leq\alpha-\beta$, and since $M_{T}$ was seen to satisfy one of the
criteria (to be precise: Criterion II) we can apply part (i) of Theorem
\ref{TCLT} and obtain%
\[
-\log\bar{F}(x)\sim px.
\]
In particular, we see that for all possible parameters in the NIG-CIR model
the formula (\ref{as}) holds true. Smile-asymptotics are now an immediate
consequence from Theorem \ref{FromBF}.

\begin{figure}[tbh]
\includegraphics{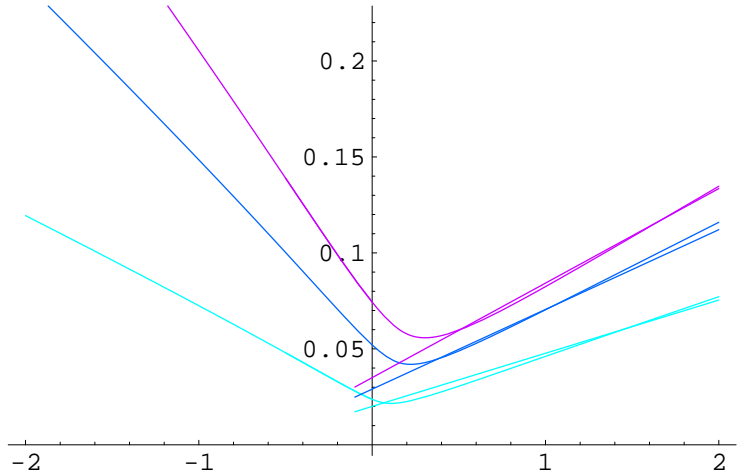}\caption{NIG with CIR time change. Parameters from
\cite{SST}. Total implied variance and slopes for three maturities $t=0.4,0.9$
and $1.3$ years. Observe that the the lines with correct slope do not
perfectly line up with the smile which is \textit{not} a contradiction to the
result that that $V^{2}\left(  k\right)  /k$ converges to a constant.}%
\end{figure}

\subsection{\bigskip The Heston Model}

The Heston model is a stochastic volatility model defined by the following
stochastic differential equations:%
\begin{align*}
\frac{dS_{t}}{S_{t}}  &  =\sqrt{v_{t}}dW_{t}^{1}\\
dv_{t}  &  =\kappa(\eta-v_{t})dt+v_{t}dW_{t}^{2}%
\end{align*}
where $d\left\langle W_{t}^{1},W_{t}^{2}\right\rangle =\rho dt$ is the
correlation of the two Brownian motions. $\log S_{t}$ therefore has the
distribution of a Brownian motion with drift $-1/2$ evaluated at a\ random
time $T\left(  \omega,t\right)  =$ $\int_{0}^{t}v_{s}ds$ with the distribution
of an integrated CIR\ process, as in the previous example. When $\rho=0$, the
L\'{e}vy process $L\equiv W^{1}$ and $T$ are independent and we can apply the
same analysis as above. Namely, the cgf of the Brownian motion with drift
speed $-1/2$ at time $1$ is%
\[
K_{L}(v)=(v^{2}-v)/2,\text{ }%
\]
so that $p_{L}=\infty,$ and $M_{T}=\exp\left(  K_{T}\right)  $ satisfies
Criterion II hence, by part (i) of Theorem \ref{TCLT},%
\[
\log\bar{F}(x)\sim-px
\]
where $p$ is determined by the equation $K_{L}\left(  p\right)  =p_{T}$. When
$\rho\leq0$, we can analyze the mgf of $\log S_{t}$ directly, and we can apply
the same reasoning as for the mgf of the CIR\ process, to deduce that
criterion II is satisfied. The distribution function for the Heston returns
hence satisfies $\log\bar{F}(x)\sim-px$ where $p$ is solution to, see
\cite{AP},%
\[
\left.  (\kappa-\rho v\theta)+(\theta^{2}(v^{2}-v)-(\kappa-\rho v\theta
)^{2})^{1/2}\cot\{(\theta^{2}(v^{2}-v)-(\kappa-\rho v\theta)^{2}%
)^{1/2}t/2\}\right\vert _{v=p}=0.
\]

When $\rho>0$, which is of little practical importance (at least
in equity markets), the mgf may explode at a different point, see
\cite{AP}, but criterion II will still be satisfied.

\begin{figure}[tbh]
\includegraphics{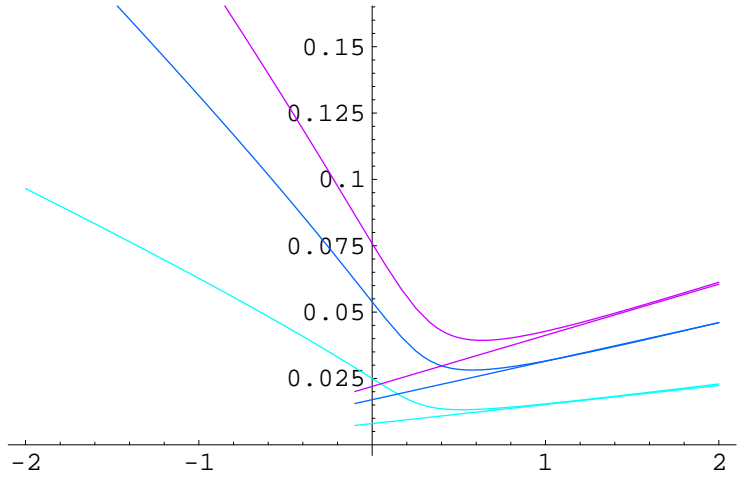}\caption{Heston Model. Parameters from \cite{SST}.
Total implied variance and slopes for three maturities $t=0.4,0.9$ and $1.3$
years.}%
\end{figure}

\begin{acknowledgement}
The authors would like to thank Chris Rogers and Nick Bingham for related
discussions. Financial support from the Cambridge Endowment for Research in
Finance (CERF) is gratefully acknowledged.
\end{acknowledgement}

\bigskip

\end{document}